# ON THE USE OF THE GENERALIZED LITTLEWOOD THEOREM CONCERNING INTEGRALS OF THE LOGARITHM OF ANALLYTICAL FUNCTIONS FOR CALCULATIONS OF INFINITE SUMS ANS ANALYSIS OF ZEROES OF ANALYTICAL FUNCTIONS

S. K. SEKATSKII

*Laboratory of Biological Electron Microscopy, IPHYS, BSP 419, Ecole Polytechnique Fédérale de Lausanne, and Dept. of Fundamental Biology, Faculty of Biology and Medicine, University of Lausanne, CH1015 Lausanne-Dorigny, Switzerland.*

*E-mail : serguei.sekatski@epfl.ch*

Abstract. Recently, we have established and used the generalized Littlewood theorem concerning contour integrals of the logarithm of analytical function to obtain a few new criteria equivalent to the Riemann hypothesis. Here, the same theorem is applied to calculate certain infinite sums and study the properties of zeroes of a few analytical functions.

## 1. INTRODUCTION

The generalized Littlewood theorem concerning contour integrals of the logarithm of analytical function is the following statement [1, 2]:

**Theorem 1.1 (The Generalized Littlewood theorem):** *Let C denote the rectangle bounded by the lines $x = X_1$, $x = X_2 > X_1$, $y = Y_1$, $y = Y_2 > Y_1$, and let f(z) be analytic and non-zero on C and meromorphic inside it. Let also g(z) be analytic on C and meromorphic inside it. Let F(z)=ln(f(z)) be the logarithm defined as follows: we start with a particular determination on $x = X_2$, and obtain the value at other points by continuous variation along y=const from $\ln(X_2 + iy)$. If, however, this path would cross a zero or pole of f(z), we take F(z) to be $F(z \pm i0)$ according as to whether we approach the path from above or below. Let also $\widetilde{F}(z) = \ln(f(z))$ be the logarithm defined by continuous variation along any smooth curve fully lying inside the contour which avoids all poles and zeroes of f(z) and starts from the same particular determination on $x = X_2$. Suppose also that the poles and zeroes of the functions f(z), g(z) do not coincide.*

*Then*

$$\int_C F(z)g(z)dz = 2\pi i \left( \sum_{\rho_g} res(g(\rho_g) \cdot \widetilde{F}(\rho_g)) - \sum_{\rho_f^0} \int_{X_1+iY_\rho^0}^{X_\rho^0+iY_\rho^0} g(z)dz + \sum_{\rho_f^{pol}} \int_{X_1+iY_\rho^{pol}}^{X_\rho^{pol}+iY_\rho^{pol}} g(z)dz \right) \quad (1.1)$$

*where the sum is over all $\rho_g$ which are poles of the function g(z) lying inside C, all $\rho_f^0 = X_\rho^0 + iY_\rho^0$ which are zeroes of the function f(z) both counted taking into account their multiplicities (that is the corresponding term is multiplied by m for a zero of the order m) and which lie inside C, and all $\rho_f^{pol} = X_\rho^{pol} + iY_\rho^{pol}$ which are poles of the function f(z) counted taking into account their multiplicities and which lie inside C. The assumption is that all relevant integrals on the right hand side of the equality exist.*

The proof of this theorem [2] is very close to the proof of the standard Littlewood theorem corresponding to the case *g(z)=1*, see e.g. [3]. In this Note we apply this Theorem for certain particular cases when the contour integral $\int_C F(z)g(z)dz$ disappears (tends to zero) if the contour tends to infinity, that is when $X_1, Y_1 \to -\infty, X_2, Y_2 \to +\infty$. This means that eq. (1.1) takes the form

$$\sum_{\rho_f^0} \int_{-\infty+iY_\rho^0}^{X_\rho^0+iY_\rho^0} g(z)dz - \sum_{\rho_f^{pol}} \int_{-\infty+iY_\rho^{pol}}^{X_\rho^{pol}+iY_\rho^{pol}} g(z)dz = \sum_{\rho_g} res(g(\rho_g)F(\rho_g)) \qquad (1.2)$$

If the integrals here can be calculated explicitly, in this way one obtains equalities involving finite or infinite sums (this last case is the most interesting one). Certainly, different methods to calculate infinite sums, and the use of these sums to study zeroes of analytical functions are well known, see e.g. [3 - 10]. Nevertheless, we believe the current paper is useful not only from the methodological and pedagogical point of view: it seems that in many cases such sums are more difficult or even impossible to calculate using other methods.

2. APPLICATION OF THE GENERALIZED LITTLEWOOD THEOREM TO CALCULATE INFINITE SUMS

Let us start with a few examples. First, we present almost trivial ones just for illustrative purposes. In a sense, the most natural function to start with is $f(z) = \cos(\pi z)$ or similar.

Evidently for large $|z|$, $\text{Im}\, z \neq 0$, $\ln(\cos(\pi z)) = O(z)$, and thus if we take $g(z) = 1/z^3$, the contour integral $\int_C \dfrac{\ln(\cos(\pi z))}{z^3} dz$ for $X_1$, $Y_1 \to -\infty$, $X_2$, $Y_2 \to +\infty$ vanishes (when constructing the contour one just should avoid such values of $X_1$, $X_2$ which give $\cos(\pm \pi X_{1,2}) = 0$). Zeroes of the function $f(z)$ lye at $z_\rho = 1/2 \pm n$ where $n$ is an integer, and $\int_{-\infty}^{1/2+n} \dfrac{1}{z^3} dz = -\dfrac{1}{2(n+1/2)^2}$. Further, the integrand has a single pole of the third order at $z=0$ and thus the corresponding residue contribution is $\dfrac{1}{2} \dfrac{d^2 \ln(\cos(\pi z))}{dz^2}\Big|_{z=0} = -\dfrac{\pi^2}{2}$. Hence the application of (1.2) gives $\sum_{n=-\infty}^{\infty} \dfrac{1}{(n+1/2)^2} = \pi^2$ which is, of course, extremely well known especially if recast as $\sum_{n=0}^{\infty} \dfrac{1}{(2n+1)^2} = \pi^2/8$.

A more general integral, viz. $\int_C \dfrac{\ln(\sin(az+b))}{(z+c)^3} dz$, where $a$, $b$, $c$ are arbitrary complex numbers, similarly give $\sum_{n=-\infty}^{\infty} \dfrac{1}{(\pi n - b + ca)^2} = \dfrac{1}{\sin^2(b-ac)}$. This is an example # 6.1.27 from [11]. Of course, here $\pi n - b + ca$ for any $n$ cannot be equal to zero. Similarly as above, when constructing the proof (demonstrating the disappearance of the contour integral), the values of $X_1$, $X_2$ corresponding to zeroes of $\sin(aX_i + b)$ should be avoided, we will not explicitly mention this any more. It is instructing to compare this above proof with much more complicated and less transparent given at four pages, pp. 219 – 222, of the famous Bromwich book [4].

Other possibility is the use of the gamma function which has simple poles at $z=0$, -1, -2,… and has no zeroes [11].

For illustration, let us consider rather general case $\prod_{n=1}^{\infty}\frac{P(n)}{Q(n)}$, that is the estimation of $\sum_{n=1}^{\infty}\ln\frac{P(n)}{Q(n)}$, where $P(z)=\sum_{i=0}^{p}a_i z^{-i}$ and $Q(z)=\sum_{i=0}^{q}b_i z^{-i}$. We require $a_0=b_0=1$ and, for convergence, $a_1=b_1$. We factorize $P(z)=a_p(-1)^p\prod_{i=1}^{p}\left(\frac{1}{r_i}-\frac{1}{z}\right)$ and $Q(z)=b_q(-1)^q\prod_{i=1}^{q}\left(\frac{1}{s_i}-\frac{1}{z}\right)$, where roots $r_i$, $s_i$ are not necessarily different. Certainly, zero is not a root of our polynomials, and we require that positive integers 1, 2, 3… also are not their roots. Trivially,

$$f(z):=\ln\left(\frac{P(z)}{Q(z)}\right)'=\sum_{i=1}^{p}\frac{r_i}{z(z-r_i)}-\sum_{i=1}^{q}\frac{s_i}{z(z-s_i)}=\sum_{i=1}^{p}\left(\frac{1}{z-r_i}-\frac{1}{z}\right)-\sum_{i=1}^{q}\left(\frac{1}{z-s_i}-\frac{1}{z}\right) \quad (2.1)$$

and we consider contour integral $\int_C f(z)\ln(\Gamma(1-z))dz$. Condition $a_1=b_1$ ensures the necessary asymptotic $O(z^{-3})$ for large $|z|$. Poles at $z=0$ of eq. (2.1) contribute nothing due to $\ln(\Gamma(1))=0$ while poles at $r_i$, $s_j$ contribute $\ln(\Gamma(1-r_i))$ and $-\ln(\Gamma(1-s_j))$. Thus we have

$$\prod_{n=1}^{\infty}\frac{P(n)}{Q(n)}=\frac{\prod_{i=1}^{q}\Gamma(1-s_i)}{\prod_{i=1}^{p}\Gamma(1-r_i)} \quad (2.2)$$

Certainly, for $P(z)=1-\dfrac{a^2}{z^2}$, $Q(z)=1$, $a \neq n$ this reduces to quite known example #89.5.11 of [5]: $\prod_{n=1}^{\infty}(1-\dfrac{a^2}{n^2}) = \dfrac{1}{\pi a}\sin(\pi a)$, which can be easily obtained by our method also starting from the contour integral $\int_C \dfrac{2a^2 \ln(\Gamma(1-z))}{z(z^2-a^2)}dz$ or $\int_C \dfrac{2a^2 \ln(\dfrac{\sin(\pi z)}{z})}{z(z^2-a^2)}dz$

(here $\dfrac{2a^2}{z(z^2-a^2)} = \dfrac{d}{dz}\ln(1-\dfrac{a^2}{z^2})$).

**Remark 2.1.** In Weisstein, Eric W. *"Infinite Product"*, from *MathWorld* A Wolfram Web Resource, https://mathworld.wolfram.com/InfiniteProduct.html we find without references $\prod_{n=1}^{\infty}\dfrac{P(n)}{Q(n)} = \dfrac{b_q \prod_{i=1}^{q}\Gamma(-s_i)}{a_p \prod_{i=1}^{p}\Gamma(-r_i)}$. This is the same as (3) due to known $(-r_i)\Gamma(-r_i) = \Gamma(1-r_i)$ [12] and conditions $a_0 = b_0 = 1$ which require $a_p(-1)^p \prod_{i=1}^{p}\dfrac{1}{r_i} = 1$ and $b_q(-1)^q \prod_{i=1}^{q}\dfrac{1}{s_i} = 1$.

Now we can consider the case when $z=k=1, 2, 3\ldots$ is a simple root of, say, $P(z)$ polynomial. We will denote it as the first root $r_1$. We begin with a contour integral $\int_C f(z)\ln\big((k-z)\cdot\Gamma(1-z)\big)dz$ so that the function under the logarithm sign is regular at $z=k$.

The l.h.s. of eq. (2) changes to $\sum_{n=1, n\neq k}^{\infty} \ln \frac{P(n)}{Q(n)}$ while at r.h.s for each pole $z=0$ in eq. (2.1) we have a contribution $-\ln k$, that is $-p \ln k$ totally; and contribution at the pole $z = r_i$ not equal to $k$ of $f(z)$ changes to $\ln((k-r_i) \cdot \Gamma(1-r_i))$. At $z=k$ the following takes place.

We know that in the vicinity of $1-k$, the residue of the gamma – function is $\frac{(-1)^{k-1}}{(k-1)!}$ hence for $z = k+\delta$ $(k-z)\Gamma(1-z) = -\delta\Gamma(1-k-\delta) \sim \frac{(-1)^{k-1}}{(k-1)!}$ and we have a contribution $\ln(\frac{(-1)^{k-1}}{(k-1)!})$. Thus total contribution of zeroes of $P(z)$ is

$$-p \ln k + \ln(\prod_{i=2}^{p}(k-r_i)\Gamma(1-r_i)) + \ln\frac{(-1)^{k-1}}{(k-1)!}$$

and that of zeroes of $Q(z)$ is

$$q \ln k - \ln(\prod_{i=1}^{p}(k-s_i)\Gamma(1-s_i)).$$

Collecting everything together, we thus formulate the following theorem.

**Theorem 2.1:** Let $P(z) = \sum_{i=0}^{p} a_i z^{-i}$ and $Q(z) = \sum_{i=0}^{q} b_i z^{-i}$, $a_0 = b_0 = 1$ and $a_1 = b_1$. Let $r_1, \ldots, r_p$ are roots of $P(z)$ and $s_1, \ldots, s_q$ are roots of $Q(z)$, not necessarily different. Let $r_1 = k$, where $k$ is a positive integer and all other roots are not equal to $k$. Then

$$\prod_{n=1, n\neq k}^{\infty} \frac{P(n)}{Q(n)} = (-1)^{p-q} k^{p-q} (k-1)! \frac{\prod_{i=1}^{q}(k-s_i)\Gamma(1-s_i)}{\prod_{i=2}^{p}(k-r_i)\Gamma(1-r_i)} \qquad (2.3)$$

For example, for $\prod_{n=1, n\neq k}^{\infty}(1-\frac{k^2}{n^2})$ we get

$$\prod_{n=1, n\neq k}^{\infty}(1-\frac{k^2}{n^2}) = \frac{(-1)^{k-1}}{2} \qquad (2.4)$$

The generalization for the case when, in addition to $r_1 = k$, some other roots are equal to integers, is straightforward.

If the series with the alternating signs are involved, the use of functions of the type $\ln(\tan(az+b))$, $\ln\frac{\Gamma(-z/2+1/2)}{\Gamma(-z/2+1)}$ might be helpful: they have alternating zeroes and poles. As an example, let us consider the series $\sum_{n=0}^{\infty}\frac{(-1)^n}{n^2a^2-b^2}$; #6.1.39 of [11]. We need to use $g(z) = \frac{d}{dz}\frac{1}{a^2z^2-b^2} = -\frac{2a^2z}{(a^2z^2-b^2)^2} = -\frac{2z}{a^2(z-b/a)^2(z+b/a)^2}$ hence we analyze the integral $\int_C \ln(\tan(\pi z/2))g(z)dz$. Its disappearance for an infinitely large contour is certain and we get, taking into account the existence of two poles of the second order of the function $g(z)$ at $z = \pm b/a$: $-\sum_{n=-\infty}^{\infty}\frac{(-1)^n}{n^2a^2-b^2} - \frac{\pi}{ab\sin(\pi b/a)} = 0$. The calculation of the values of the derivatives $\frac{d}{dz}(-\frac{2z}{a^2(z\mp b/a)^2})|_{z=\pm b/a}$ readily gives the known result $\sum_{n=0}^{\infty}\frac{(-1)^n}{n^2a^2-b^2} = -\frac{1}{2b^2} - \frac{\pi}{2ab\sin(\pi b/a)}$.

Numerous other applications of this approach to sum the infinite series can be constructed.

Other potentially useful approach is the possibility to explore the functions of more complicated arguments. For example, analyzing $\int_C \frac{\ln(\Gamma(-z^2))}{(z-c)^4} dz$ we, given the simple poles at $z = \pm\sqrt{1}, \pm\sqrt{2}...$ and double pole at $z=0$, have

$$\sum_{n=1}^{\infty}(\frac{1}{(\sqrt{n}-c)^3} + \frac{1}{(-\sqrt{n}-c)^3}) - \frac{2}{c^3} = \frac{1}{2}\frac{d^3}{dz^3}\ln(\Gamma(-z^2))|_{z=c} \text{ and thus}$$

$$\sum_{n=1}^{\infty}(\frac{1}{(\sqrt{n}-c)^3} + \frac{1}{(-\sqrt{n}-c)^3}) - \frac{2}{c^3} = 6c\psi_1(-c^2) - 4c^3\psi_2(-c^2) \quad (2.6).$$

For example, for $c=i$:

$$\sum_{n=1}^{\infty}(\frac{1}{(\sqrt{n}-i)^3} + \frac{1}{(-\sqrt{n}-i)^3}) = 2i + 6i\psi_1(1) + 4i\psi_2(1) = (2 + \pi^2 - 8\zeta(3))i \simeq 2.253i.$$

The interesting cases $c=n$ can be analyzed by putting $c = \sqrt{n} + x$ and considering the limit $x \to 0$. Let us look at an example of $c=1+x$ so that $-c^2 = -1 - 2x - x^2$. First, we, starting from the general

$$\psi(-n+x) = -\frac{1}{x} + H_n^{(1)} - \gamma + \sum_{k=1}^{\infty}(H_n^{(k+1)} + (-1)^{k+1}\zeta(k+1))x^k, \quad \text{where} \quad H_n^{(k)} = \sum_{l=1}^{n}\frac{1}{l^k} \text{ is}$$

generalized harmonic number, have:

$$\psi(-1+x) = -\frac{1}{x} + 1 - \gamma + (1+\zeta(2))x + (1-\zeta(3))x^2 + .... \text{ Whence}$$

$$\psi_1(-1+x) = \frac{1}{x^2} + 1 + \zeta(2) + 2(1-\zeta(3))x + ... \quad \text{and} \quad \psi_2(-1+x) = -\frac{2}{x^3} + 2(1-\zeta(3)) + ....$$

Further,

$$\psi_1(-c^2) = \psi(-1+(-2x-x^2)) = \frac{1}{(2x+x^2)^2} + 1 + \zeta(2) + O(x) =$$

$$\frac{1}{4x^2} - \frac{1}{4x} + \frac{19}{16} + \zeta(2) + O(x)$$

and

$$c\psi_1(-c^2) = (1+x)\psi_1(-1+(-2x-x^2)) =$$
$$(1+x)(\frac{1}{4x^2} - \frac{1}{4x} + \frac{19}{16} + \zeta(2) + O(x)) = \frac{1}{4x^2} + \frac{15}{16} + \zeta(2) + O(x). \text{ Similarly,}$$

$$\psi_2(-c^2) = \psi_2(-1+(-2x-x^2)) = \frac{2}{(2x+x^2)^3} + 2(1-\zeta(3)) + O(x)$$
$$= \frac{1}{4x^3} - \frac{3}{8x^2} + \frac{3}{8x} + \frac{27}{16} - 2\zeta(3) + O(x) \text{ and}$$

$$c^3\psi_2(-c^2) = (1+x)^3\psi_2(-1+(-2x-x^2)) =$$
$$(1+3x+3x^2+x^3)(\frac{1}{4x^3} - \frac{3}{8x^2} + \frac{3}{8x} + \frac{27}{16} - 2\zeta(3) + O(x)) =$$
$$\frac{1}{4x^3} + \frac{3}{8x^2} + \frac{31}{16} - 2\zeta(3) + O(x)$$

Application of eq. 2.6 gives

$$\sum_{n=2}^{\infty}(\frac{1}{(\sqrt{n}-1)^3} + \frac{1}{(-\sqrt{n}-1)^3}) - \frac{1}{x^3} - \frac{1}{8} - 2 = \frac{3}{2x^2} + \frac{45}{8} + \pi^2 - \frac{1}{x^3} - \frac{3}{2x^2} - \frac{31}{4} + 8\zeta(3) + O(x).$$

Finally,

$$\sum_{n=2}^{\infty}(\frac{1}{(\sqrt{n}-1)^3} + \frac{1}{(-\sqrt{n}-1)^3}) = \pi^2 + 8\zeta(3) \approx 19.486 \qquad (2.7).$$

Equations (2.6) for $c=i$ and (2.7) have been tested numerically.

The following trick also deserves to be noted. For odd functions, not only integer powers but half integer powers also may be used as arguments because single valued character of the function placed under the logarithm sign can be assured.

For this, one needs to select the function of the type $\ln\dfrac{\sin(z^{n+1/2})}{z^{n+1/2}}$ or $\ln\dfrac{\tan(z^{n+1/2})}{z^{n+1/2}}$. (We have $\dfrac{\sin(z^{n+1/2})}{z^{n+1/2}} = \sum_{k=0}^{\infty}(-1)^k \dfrac{(z^{n+1/2})^{2k+1}}{(2k+1)!\,z^{n+1/2}} = \sum_{k=0}^{\infty}(-1)^k \dfrac{z^{2kn+k}}{(2k+1)!}$, etc.).

As an example, we, starting from $\displaystyle\int_C \dfrac{1}{(z-a)^3} \ln\left(\dfrac{\sin(z^{3/2})}{z^{3/2}}\right) dz$, obtain

$$\sum_{n=1}^{\infty}\left(\dfrac{1}{((\pi n)^{2/3}-a)^2}+\dfrac{1}{((\pi n)^{2/3}e^{2\pi i/3}-a)^2}+\dfrac{1}{((\pi n)^{2/3}e^{-2\pi i/3}-a)^2}\right) = -\dfrac{3}{4\sqrt{a}}\cot(a^{3/2})+\dfrac{9}{4}\dfrac{a}{\sin^2(a^{3/2})}-\dfrac{3}{2a^2} \qquad (2.8).$$

Here $(\pi n)^{2/3}$, $(\pi n)^{2/3} e^{\pm 2\pi i/3}$ are roots of the equation $\sin(z^{3/2})=0$. Equation (2.8) has been tested numerically.

Interesting relations can be also obtained if the function $\ln(f(z)-a)$ is considered. Of course, the function $f(z)-a$ has the same poles as $f(z)$ but other zeroes $\rho_i$ of the order $k_i$ where $f(\rho_i)=a$. The most natural first example here is $\ln(e^z-a)$ where for simplicity we will take $a$ real positive. The roots of $e^z - a = 0$ are $z = \ln a \pm 2\pi n i$ where $n$ is an integer or zero, and thus we have, when analyzing $\displaystyle\int_C \dfrac{\ln(e^z-a)}{z^3} dz$: $\displaystyle\sum_{n=-\infty}^{\infty}\dfrac{1}{(\ln a + 2\pi n i)^2} - \dfrac{a}{(1-a)^2} = 0$. This is convenient to recast as

$$\sum_{n=-\infty}^{\infty}\dfrac{1}{(b+2\pi n i)^2} = \dfrac{e^b}{(1-e^b)^2} \qquad (2.9)$$

where $b$ is any real not equal to 0.

This result can be, of course, obtained "in a standard fashion" starting from

$\int_C \frac{\ln(\sinh(z/2))}{(z+b)^3} dz$. But the same idea can be interestingly used to study the solutions of

the equation $\Gamma(z) = a$ exploring the contour integral $\int_C \frac{\ln((\Gamma(z)-a) \cdot z)}{z^3} dz$, or solutions of $J_n(z) = a$ exploring the contour integral $\int_C \frac{\ln(J_n(z)-a)}{z^3} dz$, etc., see below.

From now on, we will use the following easy Lemma.

**Lemma 2.1:** *Let f(z) is an analytical function defined on the whole complex plane except possibly a countable set of points. Let also this function can be presented in some vicinity of the point z=0 by the Taylor expansion $f(z) = 1 + a_1 z + a_2 z^2 + ...$ and contour integral*

$\int_C \frac{\ln(f(z))}{z^3} dz$ *tends to zero when contour C tends to infinity (see Theorem 1.1 for the details).*

*Then for the sum over zeroes $\rho_{i,0}$ having order $k_i$ and poles $\rho_{i,pole}$ having order $l_i$ of the function*

*f(z), we have* $\sum (\frac{k_i}{\rho_{i,0}^2} - \frac{l_i}{\rho_{i,pole}^2}) = a_1^2 - 2a_2$.

*Proof:* We trivially have in some vicinity of the point z=0 the Taylor expansion

$\ln(f(z)) = a_1 z + \frac{1}{2}(-a_1^2 + 2a_2)z^2 + ...$, and now the direct application of the Theorem 1 to

the integral $\int_C \frac{\ln(f(z))}{z^3} dz$ gives the statement of the lemma.

From $\Gamma(z) = \frac{1}{z} - \gamma + (\frac{1}{2}\gamma^2 + \frac{\pi^2}{12})z + ...$ we have

$$z(\Gamma(z)-a) = 1-(\gamma+a)z+(\frac{1}{2}\gamma^2+\frac{\pi^2}{12})z^2+..., \text{ so that the application of the Lemma 2.1}$$

immediately gives $\sum \frac{k_i}{\rho_i^2} - \frac{\pi^2}{6} = (\gamma+a)^2 - \gamma^2 - \frac{\pi^2}{6}$ and $\sum \frac{k_i}{\rho_i^2} = 2a\gamma + a^2$. Here $\rho_i$ are zeroes (supposedly of the order $k_i$) of the equation $\Gamma(z)-a=0$.

**Remark 3.** We have amusing $\sum \frac{k_i}{\rho_i^2} = 0$ for the solutions $\rho_i$ of the equation $\Gamma(z) = -2\gamma$. (Remind that for a pair of the complex conjugate zeroes $\rho_{1,2} = \sigma \pm it$, one has $\frac{1}{\rho_1^2} + \frac{1}{\rho_2^2} = \frac{2(\sigma^2-t^2)}{(\sigma^2+t^2)^2}$). As follows from our analysis, this is not an isolated property but, other way around, similar properties can be established for many analytical functions. For example, as a corollary of Lemma 2.1, the equality $\sum \frac{k_i}{\rho_i^2} = 0$ holds for the solutions $\rho_i$ of the following equations: $\exp(z)-a-z-\frac{z^2}{2}=0$ ($a \neq 1$), $\exp(z)-2z=0$, $\sin(z)-a-z=0$ ($a \neq 0$), $\cos(z)-a-\frac{z^2}{2}=0$ ($a \neq 1$),

$z^{-\alpha}J_\alpha(z)-a+\frac{z^2}{4\Gamma(2+\alpha)}=0$ ($a \neq 1$, $J_\alpha(z)$ is Bessel function), and so force.

3. APPLICATION OF THE GENERALIZED LITTLEWOOD THEOREM TO ANALYZE ZEROES OF ANALYTICAL FUNCTIONS

Similar technique can be applied to calculate different sums related with zeroes of other analytical functions, that is, in a sense, not so much to study the infinite sums over the known numbers (like integers), but to analyze the properties of zeroes themselves. Numerous examples of the kind are known, see e.g. [7] and references therein. In particular, sums over non-trivial zeroes of the Riemann zeta-function, see e.g. [13] for the general discussion of this function, were extensively studied and some of them were attempted to be used to test the Riemann hypothesis [14 - 16].

By historical reasons these sums are commonly expressed via Stieltjes constants $\gamma_n$ of the Laurent expansion of the Riemann zeta-function at $z=1$: $\zeta(z) = \frac{1}{z-1} + \sum_{n=0}^{\infty} \frac{(-1)^n}{n!} \gamma_n (z-1)^n$ [13], so let us follow this line. Due to the symmetry of non-trivial Riemann zeta-function zeroes $\rho_i$, having order $k_i$, with respect to the axis Re$z=1/2$ [13], we have $\sum \frac{k_i}{\rho_i^2} = \sum \frac{k_i}{(\rho-1)_i^2}$ thus we can use the contour integral $\int_C \frac{1}{z^3} \ln(\varsigma(z-1) \cdot (z-1)) dz$ to calculate $\sum \frac{k_i}{\rho_i^2}$. Noting that

$\varsigma(z-1) \cdot (z-1) = 1 + \gamma(z-1) - \gamma_1(z-1)^2 + O((z-1))^3$, we, applying Lemma 2.1 get the well-known $\sum \frac{k_i}{\rho_i^2} = 1 - \frac{\pi^2}{8} + \gamma^2 + 2\gamma_1$, which is our first illustration. Here $\frac{\pi^2}{8} - 1 = \sum_{n=1}^{\infty} \frac{1}{(2n+1)^2}$ is the sum over the trivial zeroes of the $\zeta(z-1)$ function.

### 3.1 Incomplete gamma-function and incomplete Riemann zeta-function

As not such a trivial example, let us consider the *incomplete* gamma function, and then the *incomplete* Riemann zeta-function. Asymptotic of these functions for large arguments, see e.g. [17, 18], enables to achieve zero values of the using contour integrals used when contour tends to infinity, and we will not repeat this any more.

Incomplete gamma function is defined for Re$s$>0 and real positive $z$ as

$\gamma(s,\ z) = \int_0^z t^{s-1} e^{-t} dt$. For fixed $s$ the Taylor expansion of exponent and subsequent integration readily give well known

$$z^{-s}\gamma(s,\ z) = \sum_{n=0}^{\infty} \frac{(-1)^n}{n!(s+n)} z^n \qquad (3.1)$$

If $s \neq 0,\ -1,\ -2...$, this absolutely converged series defines the whole function on the full complex $z$-plane, and from this expression the sums $\sum_\rho \frac{k_i}{\rho_i^m}$, where $\rho_i$ is a zero of the order $k_i$, can be easily determined for all integers $m \geq 2$ using the technique of the paper. For example, from the first terms $z^{-s}\gamma(s,z) = \frac{1}{s} - \frac{z}{s+1} + \frac{z^2}{2(s+2)} + O(z^2)$, we, applying the Lemma 2.1, immediately obtain

$$\sum_{\rho_i} \frac{k_i}{\rho_i^2} = \frac{s^2}{(s+1)^2} - \frac{s}{s+2} = -\frac{s}{(s+1)^2(s+2)} \qquad (3.2)$$

This sum is never equal to zero.

Incomplete Riemann zeta-function is defined for Re$s$>1 and real positive $z$ as

$F(s,\ z) = \frac{1}{\Gamma(s)} \int_0^z \frac{t^{s-1}}{\exp(t)-1} dt$. At a moment, we regard it as a function of $z$ for some fixed $s$. We have the following series converging for $|t| < 2\pi$ and with the convention $B_1$ =-1/2; certainly, all other larger odd Bernoulli numbers are equal to zero [19]:

$$\frac{t}{e^t - 1} = \sum_{n=1}^{\infty} \frac{B_n}{n!} t^n \qquad (3.3)$$

Thus for $|z| < 2\pi$

$$z^{1-s}F(s,\ z)\Gamma(s)=\sum_{n=0}^{\infty}\frac{B_n z^n}{n!(s+n-1)} \tag{3.4}$$

here $s \neq 1,\ 0,\ -1,\ -2\ldots$ (The series representation for larger values of $|z|$ is somewhat more complicated, see e.g. [17], but this question is irrelevant for us now).

Eq. (3.1) enables to find sums $\sum_{\rho}\frac{k_i}{\rho_i^m}$, where, as usual, $\rho_i$ is a zero of incomplete gamma-function of the order $k_i$, for all integers $m \geq 2$. (Note, that already from the first terms of this development, $z^{1-s}F(s,\ z)\Gamma(s)=\frac{1}{s-1}-\frac{z}{2s}+\frac{z^2}{12(s+1)}+\ldots$ we see that for small $s$, $z$ one obtains a zero at $z \cong -2s$, and this root evidently dominates such sums for the case). From $z^{1-s}F(s,\ z)\Gamma(s)(s-1)=1-\frac{z(s-1)}{2s}+\frac{z^2(s-1)}{12(s+1)}+O(z^3)$, we get with Lemma1:

$$\sum_{\rho_i}\frac{k_i}{\rho_i^2}=\frac{(s-1)^2}{4s^2}-\frac{s-1}{6(s+1)} \tag{3.5}$$

Incomplete Riemann zeta-function at certain occasions indeed has double zeroes, see [17, 18] and references therein, hence $k_i$ are not superfluous here and below. For small $|s|$ the sum is dominated by zero at $z \simeq -2s$. It amusingly *is* equal to zero for $s_{1,2}=\pm\sqrt{3}$.

More interesting and important question is the analysis of the corresponding sums when the functions at hand are understood as functions of *s* for some fixed *z*. Eq. (3.1) immediately reveals the presence of simple poles of the incomplete gamma-function at *s*=0, -1, -2… Further, from this equation we have for $|s|<1$

$$z^{-s}\gamma(s,\ z) = \frac{1}{s} + \sum_{n=1}^{\infty} \frac{(-1)^n}{n \cdot n!} z^n (1 - \frac{s}{n} + \frac{s^2}{n^2} - ...),$$ and this absolutely converging series can be rearranged to give

$$z^{-s}s\gamma(s,\ z) = 1 + [-z + \frac{z^2}{2 \cdot 2!} - \frac{z^3}{3 \cdot 3!} + ...]s + [z - \frac{z^2}{2^2 \cdot 2!} + \frac{z^3}{3^2 \cdot 3!} - ...]s^2 + O(s^3);$$ the series, which converge for any $z$ for $|s|<1$ (their continuation for larger $s$ is irrelevant for our purposes). Thus, we should define the functions $a_m(z) := (-1)^m \sum_{n=1}^{\infty} \frac{(-1)^{n+1}}{n^m n!} z^n$ to write

$$z^{-s}s\gamma(s,\ z) = 1 + \sum_{n=1}^{\infty} a_n(z) s^n.$$ (One can add the convention $a_0(z) \equiv 1$ and write

$$z^{-s}s\gamma(s,\ z) = \sum_{n=0}^{\infty} a_n(z) s^n).$$

We obtain, using the same technique as above

$$\sum \frac{k_i}{\rho_i^2} = a_1^2(z) - 2a_2(z) + \frac{\pi^2}{6} \tag{3.6}$$

Here again the term $\frac{\pi^2}{6} = \sum_{n=1}^{\infty} \frac{1}{n^2}$ is the contribution of simple poles lying at $s$ = -1, -2, -3... For $z$ tending to zero this sum tends to $\pi^2/6 = 1.6449...$ . This reflects e.g. the known fact that for real $z<0$ incomplete gamma function has two real zeroes in each interval $-2n<s<2-2n$ for $n$=1, 2, 3… as well as known statements concerning zeroes for the case $z>0$, see [17, 18] and references therein.

**Remark 3.1:** "Almost elementary" functions can be used here. We know

$$Ein(z) := \int_0^z \frac{1-e^{-t}}{t} dt = z - \frac{z^2}{2! \cdot 2} + \frac{z^3}{3! \cdot 3} - ...,$$ whence $a_1(z) = -Ein(z)$, and further

$$\int_0^z \frac{Ein(t)}{t} dt = z - \frac{z^2}{2! \cdot 2^2} + \frac{z^3}{3! \cdot 3^2} - ... = a_2(z).$$

**Remark 3.2:** Eq. (3.3) can be exploited to determine the asymptotic of the function $a_2(z) = \int_0^z \frac{Ein(t)}{t} dt$ for large positive $z$. For such a case, we know the asymptotic $Ein(z) \sim \ln z + \gamma + o(1)$ whence $a_1^2(z) = Ein^2(z) \sim \ln^2 z + 2\gamma \ln z + \gamma^2 + o(1)$.

We have $a_2(z) = \int_0^z \frac{Ein(t)}{t} dt \sim \frac{1}{2}\ln^2 z + \gamma \ln z + O(1) + o(1)$. For large positive $z$, zeroes of the incomplete gamma function become very large in modulus and the sum $\sum \frac{k_i}{\rho_i^2}$ must tend to zero, sf. [17, 18], so that $a_1^2(z) - 2a_2(z) \sim o(1)$. Thus

$$\int_0^z \frac{Ein(t)}{t} dt \sim \frac{1}{2}\ln^2 z + \gamma \ln z + \frac{1}{2}\gamma^2 + \frac{\pi^2}{12} + o(1).$$

Analogous consideration is applicable for incomplete Riemann zeta-function. Let us now, for Re$s$>1, define the function $Q(s, z) := \Gamma(s)F(s, z) = \int_0^z \frac{t^{s-1}}{\exp(t)-1} dt$. We can rewrite (3.3) as $\frac{zt}{e^{zt}-1} = \sum_{n=0}^\infty \frac{B_n}{n!}(zt)^n$, again $B_1=-1/2$. Certainly, with a variable change $t=zy$, if Re$s>1$ and $|z|<2\pi$: $Q(s,z) = \int_0^1 (zy)^{s-1} \frac{z}{\exp(zy)-1} dy = z^{s-1}\int_0^1 y^{s-2} \frac{zy}{\exp(zy)-1} dy$, thus we obtain $z^{1-s}Q(s, z) = \sum_{n=0}^\infty \frac{B_n z^n}{n!(n+s-1)}$. For $|z|<2\pi$, this is usable for all $s$ except the simple poles at $s=1, 0, -1, -2…$. Further,

$$z^{1-s}sQ(s,\ z) = -s - s^2 - s^3 - \ldots - \frac{z}{2} + \sum_{n=2}^{\infty} \frac{B_n z^n s}{n!(n+s-1)}$$

$$= -s - s^2 - s^3 - \ldots - \frac{z}{2} + \sum_{n=2}^{\infty} \frac{B_n z^n s}{n!(n-1)}(1 - \frac{s}{n-1} + \frac{s^2}{(n-1)^2} - \ldots),\text{ so that}$$

$$-z^{1-s}sQ(s,\ z) = \frac{z}{2} + [1 - \sum_{n=2}^{\infty} \frac{B_n z^n}{n!(n-1)}]s + [1 + \sum_{n=2}^{\infty} \frac{B_n z^n}{n!(n-1)^2}]s^2 + \ldots \text{ This is valid for } |z| < 2\pi,$$

$|s| < 1$ and can be used to search for the sum of powers of inverse zeroes for such relatively small values of $z$.

For larger $z$ analytic continuation of the series over $z^n$ are required. This can be obtained similarly to the Remark above. Let us define, for $k=1, 2, 3\ldots$,

$b_k(z) := 1 + (-1)^k \sum_{n=2}^{\infty} \frac{B_n z^n}{n!(n-1)^k}$. We write $\frac{1}{t(e^t - 1)} = \sum_{n=0}^{\infty} \frac{B_n}{n!} t^{n-2}$ with $B_1 = -1/2$ and

$$\int_0^z (\frac{1}{t(e^t - 1)} - \frac{1}{t^2} + \frac{1}{2t}) dt = \sum_{n=2}^{\infty} \frac{B_n}{n!(n-1)} z^{n-1}. \text{ Thus}$$

$$b_1(z) := 1 - \sum_{n=2}^{\infty} \frac{B_n}{n!(n-1)} z^n = 1 - z \int_0^z (\frac{1}{t(e^t - 1)} - \frac{1}{t^2} + \frac{1}{2t}) dt \qquad (3.7)$$

Further, for $k=1, 2, 3\ldots$, $(-1)^k \sum_{n=2}^{\infty} \frac{B_n}{n!(n-1)^k} z^{n-2} = \frac{b_k(z) - 1}{z^2}$ and we have, by integration and multiplication on $z$, $(-1)^k \sum_{n=2}^{\infty} \frac{B_n}{n!(n-1)^k} z^n = z \int_0^z \frac{b_k(t) - 1}{t^2} dt$ whence

$$b_{k+1}(z) = 1 - \int_0^z \frac{b_k(t) - 1}{t^2} dt \qquad (3.8)$$

Now we can write

$$-z^{1-s}sQ(s,\ z) = \frac{z}{2} + b_1(z)s + b_2(z)s^2 + b_3(z)z^3 + \ldots \qquad (3.9)$$

with the aforementioned functions $b_k(z)$ defined by (3.7) and (3.8) for the whole complex plane except the points $z = 2\pi i n$ where $n$ is positive or negative integer not equal to zero. The question of the convergence of series (3.9) for $|z| >= 2\pi$ is naturally posed. We will not study this question in details: fortunately, for our current purposes it is enough that these series, for each $z$, converge in some vicinity of the point $s$=0. For this, it is sufficient that for any $z$ some positive number $A(z)$ (possibly very large), such that asymptotically for large $|z|$ $b_k$ is not larger than $A^k$, exists: than the series at question converge for $|s|<1/|A|$. This is evidently correct. For example, from (3.7) and (3.8) we have for large $z$ the asymptotic $b_1(z) \sim z \ln z$, $b_2(z) \sim \frac{1}{2} z \ln^2 z$, $b_3(z) \sim \frac{1}{6} z \ln^3 z$, etc.

This, the application of Lemma 1 to the function

$-2z^{-s}sQ(s, z) = 1 + 2z^{-1}b_1(z)s + 2z^{-1}b_2(z)s^2 + ...$ solves the problem of finding the expression for the following sum of zeroes of the function $Q(s, z)$:

$$\sum_\rho \frac{k_i}{\rho_i^2} = \frac{4}{z^2} b_1^2(z) - \frac{4}{z} b_2(z) + \frac{\pi^2}{6} + 1 \qquad (3.10)$$

Here $\frac{\pi^2}{6} + 1$ is the contribution of poles of the function $sQ(s, z)$ lying at 1, -1, -2, -3…

## 4. ZEROES OF POLYGAMMA FUNCTIONS

Recently, in interesting paper [7] Mezö and Hoffman established numerous sums including zeroes of digamma function and its generalizations. All those results can be obtained by our method, but it makes no sense to repeat. Instead, to finish the Section we concentrate on some additional possibilities brought in by the generalized Littlewood theorem.

Just to illustrate the above point, let us consider the contour integral $\int_C \frac{1}{z^3}\ln(\psi(z)\cdot(-z))dz$ with the same contour $C$ as above. We know the Laurent expansion [12, 21]

$$\psi(x) = -\frac{1}{x} - \gamma + \sum_{n=1}^{\infty}(-1)^{n-1}\zeta(n+1)x^n \qquad (4.1)$$

whence $\ln(\psi(z)\cdot(-z)) = \gamma z + (-\frac{\gamma^2}{2} - \zeta(2))z^2 + O(z^3)$. Known asymptotic of digamma function for large $z$, $\psi = O(\ln z)$ [12] guaranties that the contour integral value tends to zero when $X \to \infty$. Thus, application of the Generalized Littlewood Theorem gives:

$$0 = -\frac{\gamma^2}{2} - \zeta(2) + \frac{1}{2}\sum_{n=0}^{\infty}\frac{1}{\rho_n^2} - \frac{1}{2}\sum_{n=1}^{\infty}\frac{1}{n^2}.$$ Here the sum $-\frac{1}{2}\sum_{n=1}^{\infty}\frac{1}{n^2}$ is the contribution of simple poles of the function $\psi(x)\cdot(-x)$ at the points $z=-1, -2, -3...$ Remembering $\zeta(2) = \frac{\pi^2}{6}$, we thus have proven the relation $\sum_{i=0}^{\infty}\frac{1}{\rho_i^2} = \gamma^2 + \frac{\pi^2}{2}$ from [7]. Here we partly retain the notation of [7] and denote zeroes of digamma function, they all are real and simple, as $\rho_i$ arranged in decreasing order; $\rho_0 = 1.461632...$ is the only one positive zero.

But the aforementioned asymptotic of digamma function shows that the value of contour integral $\int_C \frac{1}{z^2}\ln(\psi(z)\cdot(-z))dz$ is also equal to zero when the contour tends to infinity. With this, the application of the Generalized Littlewood theorem, which "on equal footing" treats zeroes and poles, gives

$$\frac{1}{\rho_0} + \sum_{n=1}^{\infty}(\frac{1}{\rho_n} + \frac{1}{n}) = -\gamma \qquad (4.2)$$

Now during the summation we were obliged to group zeroes of digamma function and factors *1/n* in pairs (of course, there are many possibilities to do so). Recalling the definition $\gamma = \lim_{N \to \infty} (\ln N - \sum_{n=1}^{N} \frac{1}{n})$ [12], we can write a relation

$$\lim_{N \to \infty} (\ln N + \sum_{n=0}^{N} \frac{1}{\rho_i}) = 0 \qquad (4.3)$$

Analyzing the contour integral $\int_C \frac{1}{(z-p)^2} \ln(\psi(z)) dz$ with *p* not equal to 0, -1, -2… and not coinciding with any zero of digamma function, we obtain

$$\sum_{n=1}^{\infty} (\frac{1}{\rho_n - p} + \frac{1}{n+p}) = -\frac{\psi'(p)}{\psi(p)} \qquad (4.4)$$

In particular, using the known expressions for polygamma function on integers, for positive integer *k* we write: $\sum_{n=1}^{\infty} (\frac{1}{\rho_n - k} + \frac{1}{n+k}) = \frac{\gamma - H_{k-1}}{H_{k-1}^{(2)} - \frac{\pi^2}{6}}$, where $H_k^{(l)} = \sum_{k=1}^{k-1} \frac{1}{k^l}$ are harmonic numbers.

The case when *p* coincides with some pole or zero does not make a real problem. Let us consider, for example, $\int_C \frac{1}{(z+n)^2} \ln(-\psi(z) \cdot (z+n)) dz$. We know

$$\psi(-n+x) = -\frac{1}{x} + H_n - \gamma + \sum_{k=1}^{\infty} (H_n^{(k+1)} + (-1)^{k+1} \zeta(k+1)) x^k \text{ whence}$$

$$-x\psi(-n+x) = 1 + (\gamma - H_n)x - \sum_{k=1}^{\infty} (H_n^{(k+1)} + (-1)^{k+1} \zeta(k+1)) x^{k+1} \text{ and}$$

$$\frac{1}{\rho_k + k} + \sum_{n=1, n \neq k}^{\infty} (\frac{1}{\rho_n + k} + \frac{1}{n-k}) = H_n - \gamma \qquad (4.5)$$

In the same fashion, from $\int_C \frac{1}{(z-\rho_k)^2} \ln(\psi(z)/(z-\rho_k))dz$, we get:

$$\sum_{n=0, n\neq k}^{\infty} \left\{ \frac{1}{\rho_n - \rho_k} + \frac{1}{n+\rho_k} \right\} = -\frac{1}{2}\frac{\psi''(\rho_k)}{\psi'(\rho_k)} \qquad (4.6)$$

Similarly, for Barnes analogue of the digamma function $\psi_G(z)$, see [7]:

$$\sum_{n=1}^{\infty} \left\{ \frac{1}{\beta_n} + \frac{1}{n} \right\} = -\gamma - \frac{1}{2}\ln(2\pi) + \frac{1}{2} \qquad (4.7)$$

Here $\beta_n$ are zeroes of the Barnes digamma function: they all are real and simple [7]. For clarity of presentation we count them in simple decreasing order, so that $\beta_{1,2}$ are positive roots while others are negative. We used Lemma 3.1 from [7]: $\ln(z\psi_G(z)) = (\gamma + \frac{1}{2}\ln(2\pi) - \frac{1}{2})z + O(z^2)$ here.

Along the same line, zeroes of polygamma functions can be analyzed. For example, the tri-gamma function is defined, whenever appropriate, as $\psi'(z) = \frac{d\psi}{dz} = \sum_{n=1}^{\infty} \frac{1}{(n+z)^2}$. More generally $\psi^{(n)}(z) = \frac{d^n \psi}{dz^n} = \sum_{n=1}^{\infty} \frac{(-1)^{n+1} n!}{(n+z)^{n+1}}$. From (4.1) we have

$$\psi^{(n)}(z) = \frac{(-1)^{n+1} n!}{z^{n+1}} + \sum_{k=0}^{\infty} (-1)^{n-1+k} \frac{(n+k)!}{k!} \zeta(n+k+1) z^k \qquad (4.8)$$

Correspondingly, $\frac{1}{n!}\psi^{(n)}(z)(-1)^{n+1} z^{n+1} = 1 + \sum_{k=0}^{\infty} (-1)^k C_{n+k}^k \zeta(n+k+1) z^{k+n+1}$ where $C_{n+k}^k = \frac{(n+k)!}{n!k!}$ is a binomial coefficient.

First, we consider $\int_C \frac{1}{z^2} \ln(\psi^{(n)}(z) \cdot (-1)^{n+1} z^{n+1})) dz$. The point $z=0$ is regular but we have poles of the $n+1$ order at $z=-1, -2, -3\ldots$ Thus, we obtain the relation

$$\sum_{m=1}^{\infty} (\frac{k_{n,m}}{\eta_{n,m}} + \frac{n+1}{m}) = 0.$$ (Certainly, if we suppose that all zeroes are simple, then

$$\sum_{k=1}^{\infty} (\frac{1}{\eta_{n,k}} + \frac{n+1}{k}) = 0.$$ This, however, apparently is not proven).

Similarly, for $\int_C \frac{1}{z^{j+1}} \ln(\psi^{(n)}(z) \cdot (-1)^{n+1} z^{n+1})) dz$ with $2 \leq j \leq n$ we have:

$$\sum_{k=1}^{\infty} (\frac{l_{n,k}}{\eta_{n,k}^j} - \frac{n+1}{(-k)^j}) = 0.$$ That is $\sum_{k=1}^{\infty} \frac{l_{n,k}}{\eta_{n,r}^j} = (-1)^j (n+1) \zeta(j)$. Finally, let us consider

$\int_C \frac{1}{z^{n+2}} \ln(\psi^{(n)}(z) \cdot (-1)^{n+1} z^{n+1})) dz$. We know

$\ln(\psi^{(n)}(z)(-1)^{n+1} z^{n+1}) = \zeta(n+1) z^{n+1} + O(z^{n+2})$ so that

$$\zeta(n+1) + \frac{1}{n+1} \sum_{k=1}^{\infty} (\frac{l_{n,k}}{\eta_{n,k}^{n+1}} - \frac{n+1}{(-k)^{n+1}}) = 0.$$ For odd $n$ thus simply $\sum_{k=1}^{\infty} \frac{l_{n,k}}{\eta_{n,k}^{n+1}} = 0$. For even $n$:

$$\sum_{k=1}^{\infty} \frac{l_{n,k}}{\eta_{n,k}^{n+1}} = -2(n+1)\zeta(n+1).$$ Sums including larger powers of zeroes, like $\sum_{k=1}^{\infty} \frac{l_{n,k}}{\eta_{n,k}^{n+1+l}}$ with $l=1, 2, 3\ldots$ can be obtained in the same fashion.

For clarity of presentation, we collect most of these results as the following Theorem.

**Theorem 4.1:** *Let $\eta_{n,k}$ are zeroes of polygamma function $\psi^{(n)}(z)$ taken in increasing order of their module and $l_k$ is an order of such zero. Then*

*a)* $\sum_{k=1}^{\infty} (\frac{l_{n,k}}{\eta_{n,k}} + \frac{n+1}{k}) = 0$

b) for $2 \leq j \leq n$: $\sum_{k=1}^{\infty} \frac{l_{n,k}}{\eta_{n,r}^{j}} = (-1)^{j}(n+1)\zeta(j)$

c) For n=1, 2, 3… $\sum_{k=1}^{\infty} \frac{l_{2n+1,k}}{\eta_{2n+1,k}^{2n+2}} = 0$, $\sum_{k=1}^{\infty} \frac{l_{2n,k}}{\eta_{2n,k}^{2n+1}} = -2(2n+1)\zeta(2n+1)$.

## 5. CONCLUSIONS

We showed that the generalized Littlewood theorem concerning contour integrals of the logarithm of analytical function can be successfully applied to sum certain infinite series and analyze properties of zeroes and poles of analytical functions. In the paper, we have considered a few examples, and there is no doubt that numerous other applications of this approach will be found.

**Conflicts of Interest:** The author declares that there are no conflicts of interest regarding the publication of this paper.